# Applying GreenLab Model to Adult Chinese Pine Trees with Topology Simplification

Hong Guo[1]  Veronique Letort[2]  Xiangdong Lei[1]  Yuanchang Lu[1]  Philippe de Reffye[3]
1 Research Institute of Forest Resource Information Techniques, Chinese Academy of Forestry, Beijing, China  hongguo@caf.ac.cn, xdlei@caf.ac.cn, ylu@caf.ac.cn
2 Laboratory of Applied Mathematics Ecole Centrale Paris, Paris, France letort @mas.ecp.fr
3 CIRAD, Montpellier, France  philippe.de_reffye@cirad.fr

**Abstract**

*This paper applied the functional structural model GreenLab to adult Chinese pine trees (pinus tabulaeformis Carr.). Basic hypotheses of the model were validated such as constant allometry rules, relative sink relationships and topology simplification. To overcome the limitations raised by the complexity of tree structure for collecting experimental data, a simplified pattern of tree description was introduced and compared with the complete pattern for the computational time and the parameter accuracy. The results showed that this simplified pattern was well adapted to fit adult trees with GreenLab.*

## 1. Introduction

Functional structural tree models (FSTMs) were developed to understand the interactions and feedbacks between the morphological structure and physiological processes of individual tree growth [1]. Thus, structure design is very important in FSTMs [2].Tree architectural structure has been described at several spatial scales in different models ranging from accurate descriptions of each organ to coarse descriptions of branching systems at the individual tree level. Most of functional structural models use a very detailed representation of tree architecture based on the description of individual organs. Such detailed organ description entail complex data bases and huge computation and measurement time, which in turn leads to difficulties in applications of FSTMs to forest management.

As a generic functional structural model, GreenLab was developed to simulate plant growth at organ scale and on the basis of architectural rules and organogenesis' processes as described in details in [3, 4]. The model has been successfully applied to agricultural crops, Chinese pine saplings and beech trees [4, 5, 6, 7]. However, adult trees have much more complex architecture which makes model fitting very difficult with real topology data. The main obstacles for designing experimental protocols aiming at describing adult trees at phytomer scales are the high number of organs (several thousands) and the complexity of their topological organization. Moreover computational performance is often a limitation for any algorithm of parametric identification of tree growth models as they require running many simulations. Aiming at simplifying tree structure appropriate for different demands in forestry, the paper introduced two structural patterns and fitted them using GreenLab based on data collected on two adults Chinese pine trees (*pinus tabulaeformis* carr.). The two patterns were compared in terms of experimental work needed, of fitting results and of computational time. The simplified pattern was found adequate for adult trees in GreenLab.

## 2. Materials and Method

### 2.1. Model Structure

The GreenLab tree model was described in details in [6]. Here, the main equations and parameters are listed as follows.

GreenLab is a discrete dynamic model simulating tree growth with a time step called growth cycle (GC) and equal to one year. Carbon production $Q(i)$ during one GC $i$ is calculated as:

$$Q(i) = \frac{E(i) \cdot S_p}{r} \left(1 - \exp\left(-k \cdot \frac{S(i)}{S_p}\right)\right) \quad (1)$$

Where $r$ is defined as the hydraulic resistance of the plant to transpiration, $E(i)$ is the average potential environmental factor during GC $i$. $S_p$ is the total ground projection area available for the plant, as defined in [8], and where $k$ is the Beer–Lambert extinction coefficient related to leaf angular deviation [9]. At each GC, the carbon production is partitioned between all growing parts of the tree according to a model of sink competition. The sink of new organs $P_o(k)$ depend on their branching order ($k$) hereafter called physiological age (PA) in reference to previous works [6]. The biomass allocated to organ $o$ of PA $k$ at GC $i$ is proportional to the ratio of available biomass $Q(i-1)$ to the plant total demand $D(i)$ which is the sum of demand of new organs and of demand of rings $D_{rg}(i)$:

$$\Delta q_o(k,i) = P_o(k) \cdot \frac{Q(i-1)}{\sum_{o\in\{a,e\}} P_o(k) N_o(k,i) + D_{rg}(i)} \quad (2)$$

where $Q(i-1)$ is the total biomass production at GC $i$-1, $D_{rg}(i)$ the demand of rings at GC $i$ and $o$ represents indices for organ types (needles: $a$; internodes: $e$). $N_o(k,i)$ is the number of organs belonging to category $o$ at GC $i$ and physiological age (PA) $k$. Axes are grouped into categories according to physiological age which was related to branching order for *pinus tabuleformis* [6].

The biomass allocated to the radial growth is proportional to the ring demand $D_{rg}(i)$ which is calculated according to the ratio of biomass supply to demand as follows:

$$D_{rg}(i) = P_0 + P_1 \cdot \frac{Q(i-1)}{D(i)} \quad (3)$$

where $P_0$ is the constant sink (dimensionless) of the ring compartment, $P_1$ is the corresponding slope per mass unit (g$^{-1}$), $D(i)$ is the total demand of the plant at GC $i$. The supply to demand ratio is assumed to represent the state of trophic competition within the tree and therefore to drive the relative allocation of biomass to rings.

Biomass available for ring growth is partitioned into every living phytomer according to their position in the tree topological structure. As presented in [6], the modeling of radial growth is based on a generalisation of the Pressler law [10] that states that annual ring increment of a phytomer is proportional to the foliage quantity located above it. The biomass allocated to ring growth of a phytomer of PA $k$ and position $p$ at GC $i$ is calculated as follows:

$$q_{rg}(k,i,p) = (\frac{1-\lambda}{D_{pool}(i)} + \frac{\lambda N_a(k,i,p)}{D_{pressler}(i)}) \cdot p_{rg}(k) \cdot l(k,i) \cdot Q_{rg}(i) \quad (4)$$

$$D_{pool}(i) = \sum_{k=1}^{PA_m} N_i(k,i) \cdot p_{rg}(k) \cdot l(k,i) \quad (5)$$

$$D_{pressler}(i) = \sum_{k=1}^{PA_m} \sum_{p\in P_k} N_a(k,i,p) \cdot p_{rg}(k) \cdot l(k,i) \quad (6)$$

Where $PA_m$ is the maximal PA in the tree, $l(k,i)$ is the length of the $N_i(k,i)$ internodes of PA $k$ on plant of GC $i$ and $p_{rg}(k)$ is the linear density of sink of rings for internodes of PA $k$.

$N_a(k,i,p)$ is the number of leaves at GC $i$ above the considered internodes $P_k$ is the set of possible positions $p$ of metamers of PA $k$ in the whole plant of age $i$. The parameter $\lambda$ ($0 \leq \lambda \leq 1$) drives the relative influence of the number of leaves on the cambial growth: if $\lambda=1$, the equation (4) is similar to the Pressler lax while if $\lambda=0$, the cambial growth is uniform whatever the position of the phytomer relatively to foliage is.

## 2.2. Structural Patterns

Two structural patterns for tree topological description at different levels of details were defined in the paper with considerations of measurement feasibility and potential applications in forestry. They are pattern 1 at organ level, called substructure model and pattern 2 including compartment for all branches. The two patterns are described in figure 1.

Pattern 1: Organ level (substructure model)

In GreenLab, branches of the same PA appeared at the same GC are simulated as identical: they have same structure, biomass and geometry. Therefore the whole tree structure can be represented by a recursive equation (see [11]) and the computation of the whole tree can be highly enhanced by factorizing its elements (i.e. similar structures are computed only once). This substructure algorithm speeds up more than 1000 times of computation time needed to simulate the complex tree architecture and biomass repartition [11].

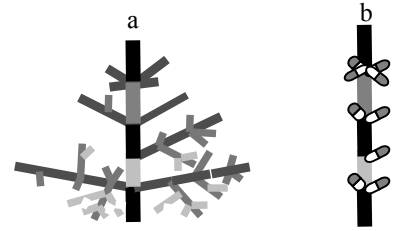

**Figure 1 Three structural patterns
(a. organ level b. compartment level )**

Pattern 2: Compartment level for branches

In pattern 2 the biomass of needles and internodes for all branches are aggregated (Figure 1b). For this pattern, stem was described internode by internode using internode length, internode radius, internode mass and needle mass, and branches were represented by total mass of needles and internodes.

### 2.3. Tree materials

Two Chinese pine trees were sampled including one 18 years old measured in autumn of 2008 and the other 31 years old measured in spring of 2006. Both were from the nursery garden located in open area in Changping forest farm, Beijing. In these experiments, the generative organs were not considered. The 3 year old needles were very scarce and most of them almost dead. So we set functioning span time for needles as 2 cycles. For the sake of simplicity, the segment of branch generated during one year was considered as a single internode although it is not true botanically speaking. At each whorl and topological level, we empirically selected one average branch which was measured internode by internode. These selected branches were described for each metamer by the following measurements: length, diameter, fresh woody weight and fresh needle weight. For other branches, we measured the total fresh weight of internodes and needles.

## 3. Data Analysis

Concerning our Chinese pine saplings, several types of input parameters listed by table 1 have been calculated directly from the measurements.

**Table 1 Directly calculated parameters and the definitions in GreenLab**

| Name | Definition |
|---|---|
| $P_e$ | Internode sink |
| $P_a$ | Blade sink |
| $b_e$ | Internode allometry |
| $\beta$ | Internode allometry |
| $\varepsilon$ | Needle SLW |

## 4. Parameters calibration

Hidden parameters showed in table 2 were estimated with nonlinear generalized least squares method [12] by Digiplante software [11].

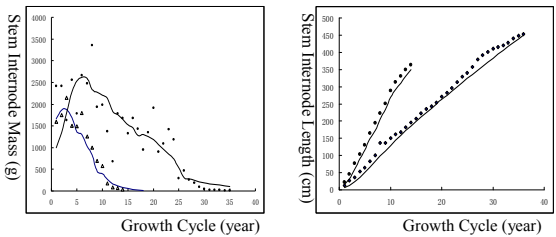

**Figure 2 Calibration of 2 trees**

Fitting results used substructure data of the two trees respectively. Figure 2 shows that the fitting results for stem internodes biomass and cumulated internodes length are consistent with the observations. Then data formatted according to pattern 2 were fitted. The computation time and error were listed in Table 3.

**Table 2 Estimated parameter values after fitting the data of Tree 1 and Tree 2**

| Name | Value | CV(%) | Value | CV(%) |
|---|---|---|---|---|
| | Tree 1 | | Tree 2 | |
| $r$ | 1.79 | 0.46 | 7.44 | 0.14 |
| $P_1$ | 0.54 | 2.6 | 0.033 | 0.30 |
| $\lambda$ | 0.01 | 6.26 | 0.40 | 0.35 |
| $p_{rg}(2)$ | 0.89 | 4.76 | 0.99 | 0.3 |
| $Sp$ | 3.04 m$^2$ | 0.49 | 78 m$^2$ | 1.0 |

**Table 3 Three patterns comparison**

| Pattern | Tree 1 | | Tree 2 | |
|---|---|---|---|---|
| | Error | Computation time (second) | Error | Computation time (second) |
| Pattern 1 | 67.3 | 1155 | 117.7 | 98780 |
| Pattern 2 | 7.2 | 180 | 6.6 | 1123 |

## 5. Results and conclusion

A potential challenge for Functional structural tree model is to face a bottleneck due to complex topological structure of adult trees [2, 14]. Substructure factorization has been introduced in detail in [12]. This paper applied substructure model and other compartment pattern to adult Chinese pine trees for comparing the results of simplification.

Firstly, based hypothesis was test for one 41 years old tree and one 18 years old tree.

We found strong adequation of the internodes allometries shape to the data for the two trees ($R^2$ ranged from 0.73 to 0.95).

Secondly, from table 2, we found that for 18 years old tree parameter $\lambda$ is equal to 0.01 which means that radial growth is nearly uniform along axes: it does not depend on the number of leaves located about each phytomer. In contrast, for 41 years old tree $\lambda=0.40$, which means that the partitioning model for radial growth is closer to the Pressler law: the relative foliage position influences for 40% of the radial growth of each phytomer. We can come to a conclusion that our general model that allows turning the relative influence of foliar biomass on ring growth is well adapted to fit Chinese pine trees of different ages.

This paper tried to find simplified topology structure to fit GreenLab model to save computation time and measurement time. Two patterns data were filled in target files and compared according to computation time and error.

From table 3, the final error of pattern 2 is lower than that of pattern 1 because the number of data is lower. Pattern 2 saved 8 times to 89 times time than Pattern 1. Thus, we can collect only compartment data to fit GreenLab model and according to different application and demand, we select different simplified way such as pattern 2. In a future work, a third pattern will be introduced at an intermediate level between pattern 1 and pattern 2. This pattern 3 would account for more detailed measurement on first order branches (internode by internode) and biomass data will be aggregated by compartment (wood/leaves) for branches of second order only. The generalization to any branching order will be considered. Thus several levels of tree description will be available in order to adapt the fitting process to flexible experimental protocols and formats of target data files.

## 6. Acknowledgements

The study was funded by Central Basic R & D Special Fund for Public Welfare Institutes of China (RIFRIGTZGZ2007008), Natural Science Foundation of China (30872022), and Chinese 948 project (2008-4-63).